\documentclass[a4paper]{article}
\usepackage[british]{babel}
\usepackage[utf8]{inputenc}
\usepackage[a4paper, total={5.5in, 8.5in}]{geometry}

\usepackage{amsmath,amssymb,amsthm}
\usepackage{dsfont}
\usepackage{tikz}
\usepackage{tikz-cd}
\usepackage[style=alphabetic,sorting=anyt,sortcites=true,backend=biber,maxbibnames=5,maxcitenames=5]{biblatex}
\usepackage{hyperref}
\usepackage{cleveref}
\usepackage{float}
\usepackage[table]{colortbl}
\usepackage{longtable}
\usepackage{csquotes}
\usepackage{caption} 

\title{\bf The Reidemeister spectra of low dimensional crystallographic groups}
\author{Karel Dekimpe\thanks{Supported by long term structural funding -- Methusalem grant of the Flemish Government. }, Tom Kaiser, Sam Tertooy\(^\ast\)}
\date{\today}

\theoremstyle{plain}
\newtheorem{theorem}{Theorem}[section]
\newtheorem{lemma}[theorem]{Lemma}
\newtheorem{proposition}[theorem]{Proposition}
\newtheorem{corollary}[theorem]{Corollary}

\theoremstyle{definition}
\newtheorem{definition}[theorem]{Definition}

\theoremstyle{remark}

\newtheorem{algorithm}{Algorithm}
\newtheorem{remark}[theorem]{Remark}

\DeclareMathOperator{\Aut}{Aut}
\DeclareMathOperator{\Out}{Out}
\DeclareMathOperator{\Inn}{Inn}

\DeclareMathOperator{\Aff}{Aff}

\DeclareMathOperator{\GL}{GL}

\DeclareMathOperator{\im}{im}
\DeclareMathOperator{\tr}{tr}

\DeclareMathOperator{\id}{id}
\newcommand{\I}{\mathds{1}}

\newcommand{\R}{\mathfrak{R}}

\DeclareMathOperator{\Spec}{Spec}

\newcommand{\NN}{\mathbb{N}}
\newcommand{\ZZ}{\mathbb{Z}}

\newcommand{\RR}{\mathbb{R}}

\newcommand{\bvarphi}{\bar{\varphi}}

\newcommand{\lie}[1]{{\mathfrak {#1}}}

\begin{document}
	
	\maketitle

	\begin{center}
	This is an Accepted Manuscript of an article published by Elsevier in Journal of Algebra on 1 Sep 2019, available online:  \href{https://doi.org/10.1016/j.jalgebra.2019.04.038}{https://doi.org/10.1016/j.jalgebra.2019.04.038}.
	
\end{center}

\begin{abstract}
	In this paper we study the number of twisted conjugacy classes (the Reidemeister number) for automorphisms of crystallographic groups. We present two main algorithms for crystallographic groups whose holonomy group has finite normaliser in \(\GL_n(\ZZ)\). The first algorithm calculates whether a group has the \(R_\infty\)-property; the second calculates the Reidemeister spectrum. We apply these algorithms to crystallographic groups up to dimension \(6\).
\end{abstract}

\section{Introduction}
In this paper we will compute Reidemeister numbers (or the number of twisted conjugacy classes) of automorphisms of crystallographic groups. Reidemeister numbers find their origin in Nielsen-Reidemeister fixed point theory \cite{jian83-1}, where they coincide with the number of fixed point classes of a self-map of a topological space.  Apart from this, Reidemeister numbers also have their importance in other fields such as Selberg theory, representation theory and algebraic geometry. We refer the reader to \cite{ft15-1} and the references therein for more information on these aspects.

The Reidemeister spectrum of a group is the collection of all Reidemeister numbers when considering  all possible automorphisms of that group. It turns out that for many groups, the only possible Reidemeister number is \(\infty\). Such groups are said to possess the \(R_\infty\)-property. This is also the case for most crystallographic groups. E.g.\ in \cite{dp11-1} it was shown that 207 out of the 219 3-dimensional crystallographic groups have the \(R_\infty\)-property. In this paper we will extend this result to the 4-dimensional groups, as well as compute the Reidemeister spectrum for groups without the \(R_\infty\)-property.

In the next section we provide the necessary preliminaries on Reidemeister numbers. Thereafter we do the same for crystallographic groups. In section 4 we present an algorithm that determines if a given crystallographic group, whose holonomy group has finite normaliser in \(\GL_n(\ZZ)\), admits the \(R_\infty\)-property. We then apply this algorithm up to dimension \(6\), and complement with some results on the remaining crystallographic groups obtained with ad-hoc methods. In section 5 we present a second algorithm, which determines the Reidemeister spectrum of a crystallographic group whose holonomy group has finite normaliser in \(\GL_n(\ZZ)\). This algorithm is then used to calculate the relevant Reidemeister spectra, up to dimension \(4\). Again, we complement with results obtained on the crystallographic groups whose holonomy groups do not have finite normaliser. Finally, we summarise all our results in the last section.

\section{Reidemeister numbers and spectrum}
In this section we introduce basic notions concerning the Reidemeister number. For a general reference on Reidemeister numbers and its connection with fixed point theory, we refer the reader to  \cite{jian83-1}.

\subsection{The Reidemeister number}
\begin{definition}
	Let \(G\) be a group and \(\varphi: G \to G\) an endomorphism. Define an equivalence relation \(\sim\) on \(G\) by
	\begin{equation*}
	\forall g,g' \in G: g \sim g' \iff \exists h \in G: g = hg'\varphi(h)^{-1}.
	\end{equation*}
The equivalence classes are called \emph{Reidemeister classes} or \emph{twisted conjugacy classes}, and we will denote the Reidemeister class of \(g\) under the endomorphism \(\varphi\) by \([g]_\varphi\). The set of Reidemeister classes of \(\varphi\) is denoted by \(\R(\varphi)\). The \emph{Reidemeister number} \(R(\varphi)\) is the cardinality of \(\R(\varphi)\) and is therefore always a positive integer or infinity.
\end{definition}

\begin{definition}
	Let \(\Aut(G)\) be the automorphism group of a group \(G\). We define the \emph{Reidemeister spectrum} as
	\begin{equation*}
	\Spec_R(G) = \{R(\varphi) \mid \varphi \in \Aut(G)\}.
	\end{equation*}
	If \(\Spec_R(G) = \{\infty\}\) we say that \(G\) has the \emph{\(R_\infty\)-property}, and if \(\Spec_R(G) = \NN \cup \{\infty\}\) we say \(G\) has \emph{full} Reidemeister spectrum.
\end{definition}

The following lemma is pivotal in determining the \(R_\infty\)-property of crystallographic groups.

\begin{lemma}[see {\cite[Lemma 1.1]{gw09-2}}]
	\label{lem:diagram3props}
	Let \(N\) be a normal subgroup of \(G\) and \(\varphi \in \Aut(G)\), with \(\varphi(N)=N\). If \(\bvarphi\) is the induced automorphism on the quotient group \(G/N\), then \(R(\varphi) \geq R(\bvarphi)\).
\end{lemma}

From this we immediately obtain:

\begin{corollary}
	\label{cor:charRoo}
	Let \(N\) be a characteristic subgroup of \(G\). If the quotient \(G/N\) has the \(R_\infty\)-property, then so does \(G\).
\end{corollary}

We will also need the following lemma, which is easy to prove.

\begin{lemma}
	\label{cor:directprod}
	Let \(G = G_1 \times G_2\) be a direct product where both \(G_1 \times \{1\}\) and \(\{1\} \times G_2\) are characteristic subgroups. Then \(\Aut(G) \cong \Aut(G_1) \times \Aut(G_2)\), and for any automorphism \(\varphi = \varphi_1 \times \varphi_2\) we have \(R(\varphi) = R(\varphi_1)\cdot R(\varphi_2)\). Hence \(\Spec_R(G) = \Spec_R(G_1)\cdot\Spec_R(G_2)\).
\end{lemma}

Finally, we give a lemma that provides equality of Reidemeister numbers of different automorphisms of the same group.

\begin{lemma}[see {\cite[Corollary 3.2]{flt08-1}}]
	\label{lem:equivreidemeisternumbers}
	Let \(G\) be a group and let \(\varphi_1,\varphi_2 \in \Aut(G)\). If there exists some \(\iota \in \Inn(G)\) such that \(\varphi_1 = \varphi_2 \circ \iota\), then \(R(\varphi_1) = R(\varphi_2)\).

\end{lemma}

\section{Crystallographic groups}
The Euclidean group \(E(n)\) is the semi-direct product \(E(n) = \RR^n \rtimes O(n)\) where multiplication is defined by \((d_1,D_1)(d_2,D_2) = (d_1+D_1d_2,D_1D_2)\). 
A cocompact discrete subgroup \(\Gamma\) of \(E(n)\) is called a \emph{crystallographic group} of dimension \(n\), or a \emph{Bieberbach group} if it is also torsion-free. 

Crystallographic groups are well understood by the three Bieberbach theorems. We refer to  \cite{wolf77-1,char86-1,szcz12-1} for more information and proofs of these theorems. The first Bieberbach theorem says that if \(\Gamma\subseteq E(n)\) is an \(n\)-dimensional crystallographic group, then its subgroup of pure translations \(T=\Gamma\cap \RR^n\) is a lattice of \(\RR^n\) (so \(T\cong \ZZ^n\)) and is of finite index in \(\Gamma\). Moreover, \(T\) is the unique maximal abelian normal subgroup of \(\Gamma\). Hence
any crystallographic group \(\Gamma\) fits in a short exact sequence
\begin{equation*}
1 \to T\cong \ZZ^n \to \Gamma \to F \to 1,
\end{equation*}
where \(T\cong\ZZ^n\) is the unique maximal abelian and normal (hence characteristic) subgroup of \(\Gamma\) and \(F\) is a finite group called the \emph{holonomy group} of \(\Gamma\). Moreover Zassenhaus showed in \cite{zass48-1} that any group \(\Gamma\) fitting in such a short exact sequence, with \(\ZZ^n\) maximal abelian in \(\Gamma\) and \(F\) finite, can be realised as a crystallographic group (i.e.\ there is an injective morphism \(i:\Gamma \to E(n)\) for which \(i(\Gamma)\) is discrete and cocompact in \(E(n)\)).

Very often, however, it is easier to realise \(\Gamma\) as a group of affine transformations. Let \(\Aff(\RR^n)=\RR^n\rtimes \GL_n(\RR)\) where multiplication is given in the same way as for \(E(n)\). Then any \(n\)-dimensional crystallographic group \(\Gamma\) can be realised as a subgroup of \(\Aff(\RR^n)\) in such a way that 
the group of pure translations \(\RR^n\cap\Gamma\) is exactly \(\ZZ^n\) (and not only isomorphic to it). It follows that any other element of \(\Gamma\) is of the form \((a,A)\) for some \(A\in \GL_n(\ZZ)\).
Let \(p:\Aff(\RR^n) \to \GL_n(\RR):(d,D) \mapsto D\) denote the projection onto the linear part, then \(F\cong p(\Gamma)\) and hence we can always view \(F\) as a subgroup of \(\GL_n(\ZZ)\). In fact \(F=\{ A\in \GL_n(\ZZ)\;|\; \exists \gamma\in \Gamma:\, \gamma=(a,A) \mbox{ for some }a \in \RR^n\}\).
In general, when we describe a crystallographic group \(\Gamma\), we will always present it by giving a collection of generators in the following form 
\begin{equation*}
\Gamma =\langle \ZZ^n, (a_1, A_1), (a_2, A_2), \ldots, (a_k,A_k) \rangle .
\end{equation*}
Here \(\ZZ^n\) stands for the subgroup of pure translations of \(\Gamma\) and is really the subgroup  of \(\Gamma\) consisting of the elements \((z,\I_n) \in \RR^n \rtimes \GL_n(\RR)\) where \(z\in \ZZ^n\) and \(\I_n\) is the \(n\times n\)--identity matrix. Recall that for the extra generators \((a_i,A_i)\) we have that \(A_i\in \GL_n(\ZZ)\).  

As we will be studying automorphisms of crystallographic groups, the following result, the second Bieberbach theorem, is crucial for us.
 
\begin{theorem}[second Bieberbach theorem]
	\label{thm:gensecbieb}
	Let \(\varphi: \Gamma \to \Gamma\) be an automorphism of a crystallographic group \(\Gamma\) with holonomy group \(F\), where \(\Gamma\subseteq \Aff(\RR^n)\) and the translation subgroup of \(\Gamma\) is exactly \(\ZZ^n\).
	 Then there exists a \((d,D) \in \Aff(\RR^n)\) such that \(\varphi(\gamma) = (d,D) \gamma (d,D)^{-1}\) for all \(\gamma \in \Gamma\). Note that \(D = \varphi|_{\ZZ^n} \) is an element of \(N_{\GL_n(\ZZ)}(F)\), the normaliser of \(F\) in \(\GL_n(\ZZ)\). 
	 
	 To shorten notation, we will write \(\varphi = \xi_{(d,D)}\) and \(N_F := N_{\GL_n(\ZZ)}(F)\).
\end{theorem}
\begin{remark}
Usually, the second Bieberbach theorem is formulated and proved in case \(\Gamma\) is a genuine subgroup of \(E(n)\). The theorem above follows easily from the usual statement, using the fact that when \(\Gamma\) is realised as a subgroup of \(\Aff(\RR^n)\), then it can be conjugated into \(E(n)\) (as any finite subgroup of \(\GL_n(\RR)\) is conjugate to a subgroup of \(O(n)\)).
\end{remark}

From now onwards we will always assume that an \(n\)-dimensional crystallographic group sits inside \(\Aff(\RR^n)\)  and that its group of pure translations is \(\ZZ^n\).

\medskip

For completeness we mention that the third Bieberbach theorem states that for any positive integer \(n\), there are only finitely many \(n\)-dimensional crystallographic groups up to isomorphism.  

\medskip

When \(\varphi\) is an automorphism of a crystallographic group, the first thing we have to check is whether or not \(R(\varphi)=\infty\). The following result answers this question completely.

\begin{theorem}[see {\cite[Corollary 3.10]{dp11-1}}]
	\label{thm:det1-AD}
	Let \(\Gamma\) be an \(n\)-dimensional crystallographic group with holonomy group \(F\) and \(\varphi = \xi_{(d,D)} \in \Aut(\Gamma)\) (where we use the notation of \cref{thm:gensecbieb}). Then
	\begin{equation*}
	R(\varphi) = \infty \iff \exists A \in F \text{ such that } \det(\I_n - AD) = 0.
	\end{equation*}
\end{theorem}
To simplify notation, we introduce the map \(|\;.\;|_\infty\) as
\begin{equation*}
|\;.\;|_\infty: \ZZ \to \NN \cup \{\infty\}: x \mapsto |x|_\infty = 
\begin{cases}
|x| & \text{ if } x \neq 0,\\
\infty & \text{ if } x = 0.
\end{cases}
\end{equation*}
In the case of torsion-free crystallographic groups, there is a handy formula to compute the Reidemeister number of an automorphism.

\begin{theorem}[averaging formula, see {\cite[Theorem 4.2]{hlp12-1}}, {\cite[Theorem 4.3]{ll09-1}}]
	\label{thm:averagingalmostbieb}
	Let \(\Gamma\) be an \(n\)-dimensional Bieberbach group with holonomy group \(F\), and \(\varphi = \xi_{(d,D)} \in \Aut(\Gamma)\). Then
	\begin{equation*}
	R(\varphi) = \frac{1}{|F|}\sum_{A \in F}|\det(\I_n - AD)|_\infty.
	\end{equation*}
\end{theorem}
In general, this formula does not hold for crystallographic groups with torsion. An extra term may be present, as we will show in \cref{prop:formulaRphi} and \cref{prop:formulaRphi2}. 

Note that the references for the above two theorems are dealing with the more general case of (coincidence) Reidemeister numbers between almost-crystallographic groups, which are subgroups of \(\Aff(G)\) for a fixed nilpotent lie group \(G\). In particular, this implies that instead of \(A\) and \(D\), one will find the symbols \(A_*\) and \(D_*\) in those references. Here the star-subscript (\(*\)) indicates the induced automorphism on the corresponding Lie algebra \(\lie{g}\) of the Lie group \(G\). In the crystallographic case however, the Lie group $G$ is the additive group \(\RR^n\) and the corresponding Lie algebra is the abelian Lie algebra \(\lie{g} = \RR^n\). This allows us to identify \(G\) and \(\lie{g}\) and then we also get that  \(A_* = A\) and \(D_* = D\).

\section{The \(R_\infty\)-property for crystallographic groups}
The \(1\)-, \(2\)- and \(3\)-dimensional crystallographic groups that do not have the \(R_\infty\)-property were determined by Dekimpe and Penninckx in \cite{dp11-1}, making use of \cite{carat06-1} to calculate the normalisers \(N_F= N_{\GL_n(\ZZ)}(F)\) and applying \cref{thm:averagingalmostbieb}. We will improve upon their algorithm and extend their results up to dimension \(6\) for crystallographic groups with finite normaliser \(N_F\), and up to dimension \(4\) for groups with infinite \(N_F\).

To create a library of crystallographic groups and calculate the normalisers \(N_F\), we used CARAT \cite{carat06-1}. Our algorithms were implemented in GAP, a system for computational discrete algebra \cite{gap18-3}, and we used the GAP-package \texttt{carat} \cite{gap18-4} to access the aforementioned library.

\medskip 

Before we set about determining which groups have the \(R_\infty\)-property, let us introduce some notation. For every \(A\) in a fixed generating set for \(F\), pick a fixed element \((a,A) \in \Gamma\). For example, we may pick the unique \(a \in \RR^n\) such that all coordinates \(a_i\) of \(a\) satisfy \(0 \leq a_i < 1\). We denote the set of these \((a,A)\) by \(F_{ext}\). Then for every element \((b,A) \in \Gamma\), we have that 
\begin{equation*}
(b,A)(a,A)^{-1} = (b-a,\I_n) \in \ZZ^n,
\end{equation*}
which means that every element of \(\Gamma\) with matrix part \(A\) is of the form \((x+a,A) = (x,\I_n)(a,A) \in \ZZ^nF_{ext}\). This also means that \(\Gamma = \langle \ZZ^n, F_{ext}\rangle\). Finally, we will use \(e_1, e_2, \dots, e_n\) to denote the elements \(((1,0,\dots, 0),\I_n)\), \(((0,1,0,\dots, 0),\I_n)\), \dots, \(((0,\dots, 0,1),\I_n)\) respectively, hence \(\ZZ^n\) (as a subgroup of \(\Gamma\)) is generated by these \(e_i\).

\medskip

For a given automorphism \(\varphi = \xi_{(d,D)}\) of a crystallographic group \(\Gamma\), we have already mentioned that \(D = \varphi|_{\ZZ_n} \in N_F\). However, the converse is not necessarily true: for a given \(D \in N_F\), there may not exist a \(d \in \RR^n\) such that the map \(\xi_{(d,D)}: \gamma \mapsto (d,D)\gamma(d,D)^{-1}\) is an automorphism of \(\Gamma\). We present an algorithm to verify whether or not such \(d\) exists, and to calculate an explicit \(d\) if it exists. The (more general) idea behind this algorithm is described in \cite[Section 4.1]{luto13-1}. 

\begin{algorithm} 
\label{alg:findd}
Given an \(n\)-dimensional crystallographic group \(\Gamma\) with holonomy group \(F\) and \(D \in N_F\). Enumerate the elements of \(F_{ext}\) by \((a_1,A_1), \ldots, (a_k,A_k)\).
\medskip

\noindent \textbf{Step 1.} The map \(F \to F: A \mapsto DAD^{-1}\) is an automorphism of \(F\). We associate a permutation \(\sigma \in \mathcal{S}_{|F|}\) to this map such that \(A_{\sigma(i)}\) is the image of \(A_i\) for every \(i\).
\medskip

\noindent\textbf{Step 2.}  If the required \(d\) were to exist, then \((d,D)(a_i,A_i)(d,D)^{-1} \in \Gamma\), or equivalently
\begin{equation}
\label{eq:existenced}
(d,D)(a_i,A_i)(d,D)^{-1}(a_{\sigma(i)},A_{\sigma(i)})^{-1} = (Da_i - a_{\sigma(i)} + (\I_n - A_{\sigma(i)})d, \I_n) \in \ZZ^n,
\end{equation} 
for all \(i = 1, \ldots, k\). Therefore, construct the matrices
\begin{equation*}
M := \begin{pmatrix}
\I_n-A_{\sigma(1)}\\
\I_n-A_{\sigma(2)}\\
\vdots\\
\I_n-A_{\sigma(k)}
\end{pmatrix} \in \ZZ^{nk \times n}, \;\; m := \begin{pmatrix}
Da_1-a_{\sigma(1)}\\
Da_2-a_{\sigma(2)}\\
\vdots\\
Da_k-a_{\sigma(k)}
\end{pmatrix} \in \ZZ^{nk},
\end{equation*}
and calculate the matrices  \(P \in \GL_{nk}(\ZZ)\), \(S \in \ZZ^{nk \times n}\) and \(Q \in \GL_{n}(\ZZ)\) such that \(S\) is the Smith normal form of \(M\) and \(PMQ = S\).
With these matrices known, calculate \(t := Pm\) and define \(d': = Q^{-1}d\), and observe that condition \eqref{eq:existenced} is equivalent to
\begin{equation}
\label{eq:existencedprime}
t + Sd' \in \ZZ^{nk}.
\end{equation}
Let \(s_1, s_2, \ldots, s_r\) be the (non-zero) invariant factors of \(S\). Since these are unique up to sign, we may assume that they are positive. Writing out the coordinates of \(t + Sd'\), we find that condition \eqref{eq:existencedprime} means that \(t_i + s_id'_i \in \ZZ\) for \(i = 1, \ldots, r\) and \(t_i \in \ZZ\) for \(i = r+1, \ldots, nk\).
\medskip

\noindent \textbf{Step 3.} We verify if \(t_{r+1}, \ldots, t_{nk} \in \ZZ\). If this is not the case, the required \(d\) does not exist and we go to \textbf{Stop}. Otherwise, we set \(d'_i = -t_i/s_i\) for \(i = 1, \ldots, r\) and \(d'_i = 0\) for \(i = {r+1}, \ldots, n\), and calculate \(d = Qd'\). The map \(\xi_{(d,D)}: \Gamma \to \Gamma: \gamma \mapsto (d,D)\gamma(d,D)^{-1}\) is then an automorphism.
\medskip

\noindent \textbf{Stop.}
\end{algorithm}

\bigskip 

The following algorithm allows us to determine whether a crystallographic group with finite normaliser \(N_F\) has the \(R_\infty\)-property or not. This is basically the algorithm from \cite{dp11-1} combined with \cref{alg:findd}, meaning no work has to be done by hand anymore.

\begin{algorithm}\label{alg:mainalg1}  Given an \(n\)-dimensional crystallographic group \(\Gamma\) with holonomy group \(F\).
\medskip
	
\noindent\textbf{Step 1.} Retrieve the normaliser \(N_F\) and check if it is finite. If it is finite, proceed to step 2, otherwise go to \textbf{Stop}.
\medskip

\noindent\textbf{Step 2.} Check for every \(D \in N_F\) whether
\begin{enumerate}
	\item there exists a \(d \in \RR^n\) such that the map \(\xi_{(d,D)}: \gamma \mapsto (d,D)\gamma(d,D)^{-1}\) is an automorphism, using \cref{alg:findd}, and
	\item \(\det(\I_n - AD) \neq 0\) for all \(A \in F\).
\end{enumerate}
The algorithm stops once a suitable \(D\) has been found, or when we have tried all \(D \in N_F\). By \cref{thm:det1-AD}, \(\Gamma\) has the \(R_\infty\)-property if and only if no such \(D\) exists.
\medskip

\noindent\textbf{Stop.}

\medskip

\end{algorithm}

We have applied \cref{alg:mainalg1} for all crystallographic groups up to dimension \(6\); the results can be found in \cref{tbl:alg2results}. Here one sees that the majority of crystallographic groups have finite normaliser \(N_F\). However, this algorithm fails when \(N_F\) is infinite, in which case we can try two things:
\begin{enumerate}
	\item Show that \(\Gamma\) does not have the \(R_\infty\)-property using \cref{thm:det1-AD}, by finding a suitable \(D \in N_F\) as above. Practically, we calculate a generating set of \(N_F\) and try out \(D\) which are words in these generators with short length.
	\item Show that \(\Gamma\) has the \(R_\infty\)-property using \cref{cor:charRoo}, by finding a characteristic subgroup \(\Gamma'\) such that \(\Gamma/\Gamma'\) is a crystallographic group (of lower dimension) with the \(R_\infty\)-property.
\end{enumerate}

We did this for the crystallographic groups with infinite normaliser up to dimension \(4\), and always found the required \(D\) or \(\Gamma'\). Therefore, we have completely determined which crystallographic groups up to dimension \(4\) have the \(R_\infty\)-property. The complete results can be found in \cref{tbl:roodim3} for dimensions \(1\) to \(3\) and in \cref{tbl:roodim4} for dimension \(4\). In these tables, ``BBNWZ'' stands for the classification system used in \cite{bbnwz78-1}, ``CARAT'' stands for the classification system used by CARAT \cite{carat06-1} and ``IT'' stands for the classification system used in the International Tables in Crystallography \cite{aroy16-1} (only up to dimension \(3\)).

\begin{table}[h]
	\centering
	\begin{tabular}{l|l|l|l}
		dim & \# groups & \# \(|N_F| < \infty\)	& \# \(|N_F| < \infty\)\\
			&			&						& \& no \(R_\infty\)\\
		\hline
		&&& \\[\dimexpr-\normalbaselineskip+2pt]
		1 & 2			& 2			& 1 	\\
		2 & 17			& 15		& 1		\\
		3 & 219			& 204		& 7		\\
		4 & 4783		& 4388		& 45	\\
		5 & 222018		& 204768	& 146	\\
		6 & 28927915 	& 26975265	& 321
	\end{tabular}
	\caption{Results of \cref{alg:mainalg1}}
	\label{tbl:alg2results}
\end{table}

\section{Calculation of Reidemeister spectra}
In this section, we will calculate the Reidemeister spectra of crystallographic groups that do not possess the \(R_\infty\)-property. 

For groups with finite \(N_F\), we will first present an algorithm to calculate \(R(\varphi)\) for a given automorphism \(\varphi = \xi_{(d,D)}\). Second, we will show in \cref{thm:finitenrReidnrs} that for every \(D \in N_F\), the set \(\{R(\varphi) \mid \varphi \in \Aut(\Gamma) \textrm{ with }\varphi|_{\ZZ^n} = D\}\) is finite, so that we may apply this algorithm to a finite number of automorphisms to find the Reidemeister spectrum.

If the normaliser is infinite, we will have to proceed by hand, which is feasible up to dimension \(3\). For dimension \(4\), we limit ourselves to calculating the spectra of only a small number of groups, mainly those where we can apply \cref{cor:directprod} or \cref{thm:averagingalmostbieb}.

\subsection{Groups with finite \(N_F\)}
We first introduce a well-known lemma, which is easy to prove (for example by using the Smith normal form):
\begin{lemma}
	\label{lem:AutZnisdet}
	Let \(A \in \ZZ^{n \times n}\). Then the number of cosets of \(\im(A)\) in \(\ZZ^n\) is \(|\det(A)|_\infty\).
\end{lemma}
The following algorithm determines the Reidemeister number \(R(\varphi)\) of a given automorphism \(\varphi = \xi_{(d,D)}\).

\begin{algorithm}\label{alg:Reidnr}
	Given an \(n\)-dimensional crystallographic group \(\Gamma\) with holonomy group \(F\) and \(\varphi = \xi_{(d,D)} \in \Aut(\Gamma)\).
\medskip

\noindent\textbf{Step 1.} Calculate \(\det(\I_n - AD)\) for all \(A \in F\). If one of these determinants is zero, then \(R(\varphi) = \infty\) by \cref{thm:det1-AD} so we go to \textbf{Stop}. If not, proceed.
\medskip

\noindent\textbf{Step 2.} For each \(A \in F\), choose a set of representatives \(x\) of the cosets \(\ZZ^n / \im(\I_n - AD)\). Then every element \((y+a,A) \in \Gamma\) will be Reidemeister equivalent to \((x+a,A)\), where \(x\) is the representative of the coset \(y + \im(\I_n - AD)\). This follows from
\begin{align*}
x - y \in  \im(\I_n - AD) &\iff \exists z \in \ZZ^n: x-y = (\I_n - AD)z\\
& \iff \exists z \in \ZZ^n: (x+a,A) = (z,\I_n)(y+a,A)\varphi(z,\I_n)^{-1}.
\end{align*}
By \cref{lem:AutZnisdet}, the number of cosets is \(|\det(\I_n - AD)|\), and \(F\) is a finite group, hence
\begin{equation}
\label{eq:cosetreps}
\bigcup_{A \in F} \{(x+a,A) \mid x + \im(\I_n - AD) \in \ZZ^n / \im(\I_n - AD)\}
\end{equation}
is a finite set.
\medskip

\noindent\textbf{Step 3.} Reduce this set to representatives of distinct Reidemeister classes. This can be done by checking pairwise if two elements are Reidemeister equivalent, and removing one of the elements from the collection if they are. 

The pairwise check can be done as follows. Let \((x+a,A)\) and \((y+b,B)\) be two elements of the set \eqref{eq:cosetreps}. They are equivalent if and only if there exists some \((z+c,C) \in \Gamma\) such that
\begin{equation*}
(x+a,A) = (z+c,C)(y+b,B)(d,D)(z+c,C)^{-1}(d,D)^{-1}.
\end{equation*}
Working out this product and considering both components separately, we arrive at the following two conditions:
\begin{enumerate}
	\item \(A = CBDC^{-1}D^{-1}\),
	\item \(x+a = (\I_n - CBDC^{-1})(z+c) + C(y+b) + (CB-CBDC^{-1}D^{-1})d\).
\end{enumerate}
Practically, we will start by checking if there exists some \(C \in F\) satisfying the first condition. We then use the fact that the first condition is satisfied to simplify the statement of the second condition, which then becomes
\begin{equation*}
\exists z \in \ZZ^n: (\I_n - AD)^{-1}\left(x+a-C(y+b) - (CB-A)d \right) - c = z,
\end{equation*}
and this can easily be verified by a computer.
\medskip

\noindent\textbf{Step 4.} Count the number of elements left in the set after the pairwise elimination, which is then the Reidemeister number \(R(\varphi)\).
\medskip

\noindent\textbf{Stop.} 
\end{algorithm}
\bigskip

While we can now calculate the Reidemeister number for a given automorphism \(\varphi = \xi_{(d,D)}\), we only know the normaliser \(N_F\) but not the automorphism group \(\Aut(\Gamma)\). As such, for a given \(D \in N_F\), we need to be able to calculate all automorphisms \(\varphi\) such that \(\varphi|_{\ZZ^n} = D\). One approach to this is to determine all \(d \in \RR^n\) such that the map \(\xi_{(d,D)}: \Gamma \to \Gamma: \gamma \mapsto (d,D)\gamma(d,D)^{-1}\) is an automorphism. 

We observe that if we find two such \(d\)'s, say \(d_1\) and \(d_2\), then \(\xi_{(d_1,D)} \circ \xi_{(d_2,D)}^{-1} = \xi_{(d_1-d_2,\I_n)}\) is an automorphism as well. Hence it suffices to find all \(d\) such that \(\xi_{(d,\I_n)}\) is an automorphism, and for every \(D \in N_F\) a single \(d\) such that \(\xi_{(d,D)}\) is an automorphism. The following theorem will show how we can do the former, while the latter can be done with \cref{alg:findd}.

\begin{theorem}
	\label{thm:dbase}
	There exists a finite set of elements \(\Delta^{base} \subseteq  \RR^n\) such that every \(\varphi \in \Aut(\Gamma)\) with \(\varphi|_{\ZZ^n} = \I_n\) is the composition of some inner automorphism \(\iota = \xi_{(d^{int},\I_n)}\) with some automorphism \(\xi_{(d^{base},\I_n)}\), where \(d^{int} \in \ZZ^n\) and \(d^{base} \in \Delta^{base}\).
\end{theorem}
This theorem can be interpreted as a description of the cohomology group \(H^1(F,\ZZ^n)\), a detailed treatise can be found in \cite[Section 4.2]{luto13-1}. We present our own proof below, since it is used to introduce some notation needed later in this section.   
\begin{proof}
	We start by finding a necessary and sufficient condition on \(d\) for \(\xi_{(d,\I_n)}\) to be an automorphism. Following the steps and notation from \cref{alg:findd}, we obtain the condition \(Sd' \in \ZZ^n\), which translates to \(d_i' \in \frac{1}{s_i}\ZZ\) for every \(i = 1, \ldots, r\).
	
	Now, take the set of all \(d'\) with \(d'_i \in \{0, \frac{1}{s_i}, \frac{2}{s_i}, \ldots, \frac{s_i-1}{s_i}\}\) for \(i = 1, \ldots, r\), and \(d_i' = 0\) for \(i = r+1, \ldots, n\). Multiplying every element of this set on the left with \(Q\), we obtain the desired set \(\Delta^{base}\).
	
	To prove this truly is the desired set, pick any \(d\) such that \(\xi_{(d,\I_n)}\) is an automorphism. Let \(d' = Q^{-1}d\), then necessarily \(d'_i \in \frac{1}{s_i}\ZZ\) for \(i = 1, \ldots, r\). Now we decompose \(d'\) as
	\begin{equation*}
	\begin{pmatrix}
	d'_1\\
	\vdots\\
	d'_r\\
	d'_{r+1}\\
	\vdots\\
	d'_{n}
	\end{pmatrix} = \begin{pmatrix}
	d'_1 - \lfloor d'_1 \rfloor\\
	\vdots\\
	d'_r - \lfloor d'_r \rfloor\\
	0\\
	\vdots\\
	0
	\end{pmatrix} +  \begin{pmatrix}
	\lfloor d'_1 \rfloor\\
	\vdots\\
	\lfloor d'_r \rfloor\\
	0\\
	\vdots\\
	0
	\end{pmatrix} + \begin{pmatrix}
	0\\
	\vdots\\
	0\\
	d'_{r+1}\\
	\vdots\\
	d'_{n}
	\end{pmatrix}.
	\end{equation*}
	By definition, the first vector must be \(Q^{-1}d^{base}\) for some \(d^{base} \in \Delta^{base}\), the second vector (call it \(d'\vphantom{x}^{int}\)) is an element of \(\ZZ^n\), and the third vector (call it \(d'\vphantom{x}^{rem}\)) satisfies \(Sd'\vphantom{x}^{rem} = 0\). Setting \(d^{int} = Qd'\vphantom{x}^{int}\) and \(d^{rem} = Qd'\vphantom{x}^{rem}\), we find that \(d^{int} \in \ZZ^n\) and \((\I_n - A)d^{rem} = 0\) for all \(A \in F\). This latter observation means that \((d^{rem},\I_n)(a,A)(d^{rem},\I_n)^{-1} = (a,A)\) for all \((a,A) \in \Gamma\).
	 Thus, we may conclude that
	\begin{align*}
	\xi_{(d,\I_n)} &= \xi_{(d^{rem},\I_n)} \circ \xi_{(d^{int},\I_n)} \circ \xi_{(d^{base},\I_n)}\\
	&= \id \circ \;\iota \circ \xi_{(d^{base},\I_n)},
	\end{align*}
	which is the required decomposition.	
\end{proof}

\begin{corollary}
	\label{cor:decompxi}
	Let \(\varphi, \psi \in \Aut(\Gamma)\) with \(\varphi|_{\ZZ^n} = \psi|_{\ZZ^n}\). Then there exists a \(d^{base} \in \Delta^{base}\) and an inner automorphism \(\iota\) such that
	\begin{equation*}
	\psi = \iota \circ \xi_{(d^{base},\I_n)} \circ \varphi.
	\end{equation*}
\end{corollary}
This decomposition means that we only need to consider finitely many \(d \in \RR^n\) for every \(D \in N_F\). The following theorem formalises this. 
\begin{theorem}
	\label{thm:finitenrReidnrs}
	Let \(D \in N_F\). Then the set \(\{R(\varphi) \mid \varphi \in \Aut(\Gamma) \textrm{ with }\varphi|_{\ZZ^n} = D\}\) is finite.
\end{theorem}
\begin{proof}
	If no \(d \in \RR^n\) exists such that \(\xi_{(d,D)}\) is an automorphism, the set is empty. Otherwise, fix such \(d\) and let \(\varphi\) be any automorphism with \(\varphi|_{\ZZ^n} = D\). Using \cref{cor:decompxi} we know that there exist \(d^{base} \in \Delta^{base}\) and \(\iota \in \Inn(\Gamma)\) such that
	\begin{equation*}
	R(\varphi) = R(\iota \circ \xi_{(d^{base},\I_n)} \circ \xi_{(d,D)}) = R(\xi_{(d^{base},\I_n)} \circ \xi_{(d,D)}),
	\end{equation*}
	where the last equality is given by \cref{lem:equivreidemeisternumbers}. Thus
	\begin{equation*}
	\{R(\varphi) \mid  \varphi \in \Aut(\Gamma) \textrm{ with }\varphi|_{\ZZ^n} = D\} = \{R(\xi_{(d^{base},\I_n)} \circ \xi_{(d,D)}) \mid d^{base} \in \Delta^{base}\},
	\end{equation*}
	which is finite because \(\Delta^{base}\) is finite.
\end{proof}

This allows us to calculate the entire Reidemeister spectrum of a crystallographic group \(\Gamma\) with finite normaliser \(N_F\). First, calculate the set \(\Delta^{base}\) as in the proof of \cref{thm:dbase}. Then, for every \(D \in N_F\), run \cref{alg:findd} to find a \(d\) such that \(\xi_{(d,D)}\) is an automorphism; and (if such \(d\) was found) use \cref{alg:Reidnr} to calculate \(R(\xi_{(d^{base}+d,D)})\) for every \(d^{base} \in \Delta^{base}\). A GAP-implementation of \cref{alg:Reidnr} produced the results found in tables \ref{tbl:roodim3} and \ref{tbl:roodim4}. Note that we have omitted the value \(\infty\) from the Reidemeister spectra in these tables. The Bieberbach groups are indicated by a star (\(*\)). 

From \cref{thm:finitenrReidnrs} we may also conclude the following:
\begin{theorem}
	Let \(\Gamma\) be a crystallographic group with holonomy group \(F\). If \(N_F\) is finite, then the Reidemeister spectrum \(\Spec_R(\Gamma)\) is finite as well.
\end{theorem}
Note that the above theorem follows more directly from the fact that \(N_F\) is finite if and only if \(\Out(\Gamma)\) is finite (see \cite[page 69]{szcz12-1}) combined with \cref{lem:equivreidemeisternumbers}.

\subsection{Groups with infinite \(N_F\)}
The Reidemeister spectra of crystallographic groups with infinite \(N_F\) have to be calculated by hand. We will start off with two families: \(\ZZ^n\) and \(\langle \ZZ^n, (0,-\I_n)\rangle\) for \(n \geq 2\). Both of these families have \(N_F = \GL_n(\ZZ)\). In fact, \(\Aut(\ZZ^n) = \GL_{n}(\ZZ)\) and \(\Aut(\langle \ZZ^n, (0,-\I_n)\rangle) = \{ \xi_{(d,D)} \mid d \in (\frac{1}{2}\ZZ)^n, D \in \GL_n(\ZZ)\}\).

\begin{proposition}[see \cite{roma11-1}]
	\label{prop:spectrumZZn}
	Let \(\Gamma = \ZZ^n\) with \(n \geq 2\). Then \(\Spec_R(\Gamma) = \NN \cup \{\infty\}\).
\end{proposition}

\begin{proposition}
	\label{prop:spectrumZZn-I}
	Let \(\Gamma = \langle \ZZ^n, (0,-\I_n)\rangle\) with \(n \geq 2\). Then
	\begin{equation*}
	\Spec_R(\Gamma) = \begin{cases}
	2\NN \cup \{3,\infty\} & \text{ if } n = 2,\\
	\NN \setminus \{1\} \cup \{\infty\} & \text{ if } n \geq 3.
	\end{cases}
	\end{equation*}
\end{proposition}

The proof of this proposition is far from straightforward. We will first introduce some lemmas and intermediate results.

\begin{lemma}
	\label{lem:OAaoddeven}
	Let \(B\in \ZZ^{n \times n}\) and \(b\in \ZZ^n\).
	Define \(O(B,b)\) as the number of solutions \(\bar{x}\) over \(\ZZ_2\) of the linear system of equations \(\bar{B}\bar{x} = \bar{b}\), where the bar-notation stands for the element-wise projection to \(\ZZ_2\). Then we have the following:
	\begin{itemize}
		\item when \(\det(B)\) is odd, \(O(B,b) = 1\) (so is also odd),
		\item when \(\det(B)\) is even, \(O(B,b) = 0,2,4, \dots, 2^n\) (so is also even).
	\end{itemize}
\end{lemma}

\begin{lemma}
	\label{lem:equivgamma2}
	Define an equivalence relation on \(\ZZ^n\), determined by a matrix \(B \in \ZZ^{n \times n}\) and an element \(b \in \ZZ^n\), where
	\begin{equation*}
	\forall x,y \in \ZZ^n: x \sim y \iff \exists z \in \ZZ^n: x-y = Bz \text{ or } x+y+b = Bz.
	\end{equation*}
	The number of equivalence classes is then given by
	\begin{equation*}
	E(B,b) = \	\frac{|\det(B)|_\infty + O(B,b)}{2}.
	\end{equation*}
\end{lemma}
\begin{proof}
It is obvious from the definition that 
\begin{equation*}
x \sim y \iff y \in x + \im(B) \mbox{ or } y \in -x-b + \im(B).
\end{equation*}
From this it follows easily that the equivalence class of \(x\), denoted by \([x]_\sim\), equals
\begin{equation}\label{union}
 [x]_\sim =  (x + \im(B)) \cup ( -x-b + \im(B)).
\end{equation}
Moreover we have that either \((x + \im(B)) \cap  ( -x-b + \im(B)) =\emptyset
 \) or \(x + \im(B)=  -x-b + \im(B)\). From \cref{lem:AutZnisdet} we know there are \(|\det(B)|_\infty\) cosets of \(\im(B)\). 

In general, elements \(x\) and \(-x-b\) will not belong to the same coset of \(\im(B)\) and the union in \eqref{union} will be a disjoint union.	
Let \(N\) denote the number of cosets \(x+\im(B)\) such that \(x+\im(B)=-x-b+\im(B)\). Then these \(N\) cosets form \(N\) equivalence classes for the relation \(\sim\), while the other \(|\det(B)|_\infty-N\) cosets come in pairs \((x + \im(B),-x-b+\im(B))\) and so determine the remaining \((|\det(B)|_\infty-N)/2\) equivalence classes of \(\sim\). Therefore
\begin{equation*}
E(B,b)= \frac{|\det(B)|_\infty-N}{2} +N =  \frac{|\det(B)|_\infty+N}{2}.
\end{equation*}
	
We now determine this number \(N\). We have that \(x + \im(B)\) and \(-x-b + \im(B)\) are actually the same coset if and only if 
\begin{equation}
\label{eq:gamma2ding}
\exists z \in \ZZ^n: 2x+b = Bz.
\end{equation}
		
We have to count for how many cosets \(x + \im(B)\) this equation holds. For it to hold, it must definitely do so over \(\ZZ_2\), i.e. \(\bar{B}\bar{z} = \bar{b}\). So we have \(O(B,b)\) solutions \(\bar{z}\) over \(\ZZ_2\). Next, we show that each solution \(\bar{z}\) over \(\ZZ_2\) produces a unique coset \(x+\im(B)\)  satisfying equation \eqref{eq:gamma2ding}. Let \(\bar{z}\) be a solution of \(\bar{B}\bar{z} = \bar{b}\). Choose any lift \(z\in \ZZ^n\) of \(\bar{z}\), 
then \(Bz -b \in 2 \ZZ^n\), so there exists a unique \(x\in\ZZ^n\) such that \(Bz -b =2x\). Hence for this \(x\) we have that equation \eqref{eq:gamma2ding} holds and so \(x+\im(B) = -x -b + \im(B)\). However, the \(x\) we found depends on the choice of the lift \(z\). Let \(z'\in \ZZ^n\) be another element projecting down to \(\bar{z}\) (so  there exists a \(c\in \ZZ^n\) with \(z-z'=2 c\)) and giving rise to \(x'\) satisfying \(2 x'= Bz' +b\). Then
\begin{equation*}
2 (x-x') = Bz -b - (Bz'-b) = 2 B c \implies x-x' = Bc \implies x + \im(B) = x' + \im(B),
\end{equation*}
from which we see that the choice of the lift \(z\) is of no influence on the coset \(x+\im(B)\): while \(x\) and \(x'\) may be different, they are both representatives of one and the same coset.

 Hence every solution \(\bar{z}\) gives rise to a unique coset with representative \(x\) satisfying equation \eqref{eq:gamma2ding}. Note that if \(\det(B) \neq 0\), each solution \(\bar{z}\) produces a different coset: suppose by contradiction that two different solutions \(\bar{z}_1\) and \(\bar{z}_2\) produce the same coset \(x + \im(B)\). This means there exist \(z_1, z_2 \in \ZZ^n\), with \(z_1\neq z_2\), such that \(2x+b = Bz_1 = Bz_2\), but then \(B(z_1-z_2) = 0\) and therefore \(\det(B) = 0\), which we assumed was not the case. So the number \(N\) of cosets \(x + \im(B)\) satisfying equation \eqref{eq:gamma2ding} is exactly \(O(B,b)\) when \(\det(B) \neq 0\). So in case   \(\det(B) \neq 0\), we have that 
\begin{equation*}
E(B,b) =  \frac{|\det(B)|_\infty+O(B,b)}{2}.
\end{equation*}
If \(\det(B) = 0\), there are infinitely many cosets \(x+\im(B)\), so there are infinitely many pairs of disjoint cosets that together form one equivalence class, and at most \(O(B,b)\) cosets that form an equivalence class on their own. Hence \(E(B,b)\) is infinite and the formula above also holds in this case.
\end{proof}

\begin{proposition}
	\label{prop:formulaRphi}
	Let \(\Gamma = \langle \ZZ^n, (0,-\I_n)\rangle\) for \(n \geq 2\), and \(\varphi = \xi_{(d,D)}\in \Aut(\Gamma)\). Then the Reidemeister number of \(\varphi\) is given by
	\begin{equation}
	\label{eq:Rgamma2}
	R(\varphi) = \left(\frac{1}{2} \sum_{A \in F} |\det(\I_n-AD)|_\infty\right) + O(\I_n-D,2d).
	\end{equation}
\end{proposition}
	
\begin{proof}

The holonomy group of \(\Gamma\) is given by \(F = \{\pm \I_n\} = Z(\GL_n(\ZZ))\). Let \(\varphi  = \xi_{(d,D)}\) be an automorphism. Note that necessarily \(d \in \left(\frac{1}{2}\ZZ\right)^n\), whereas \(D\) can be any matrix in \(\GL_n(\ZZ)\). Two elements \((x,A_x),(y,A_y) \in \Gamma \) are Reidemeister equivalent if and only if there exists an element \((z,A_z) \in \Gamma\) such that	
\begin{align*}
(y,A_y) 
&= (z,A_z)(x,A_x)\varphi(z,A_z)^{-1}\\
&= (z,A_z)(x,A_x)(d,D)(z,A_z)^{-1}(d,D)^{-1}\\
&= (z+A_zx+A_zA_xd-A_xDz-A_xd,A_zA_xDA_z^{-1}D^{-1})\\
&= (A_zx + (\I_n-A_xD)z -(\I_n-A_z)(A_xd),A_x).
\end{align*}
Thus a necessary requirement for \((x,A_x)\) to be equivalent to \((y,A_y)\) is that \(A_x = A_y\). So an element \((x,\I_n)\) and an element \((y,-\I_n)\) can never be in the same Reidemeister class, and in particular \(R(\varphi) \geq 2\). Now for two elements \((x,A), (y,A)\) with the same holonomy part \(A\), \((x,A) \sim (y,A)\) if and only if there exists some \(z \in \ZZ^n\) such that
\begin{equation*}
x-y = (\I_n-AD)z \text{ or }  x+y+2Ad = (\I_n-AD)z,
\end{equation*}
where the first case corresponds to \(A_z = \I_n\) and the second case to \(A_z = -\I_n\). From the definition of \(E(B,b)\) in \cref{lem:equivgamma2}, we obtain that
\begin{align*}
R(\varphi) &= E(\I_n-D,2d) + E(\I_n+D,-2d)\\
&= \frac{|\det(\I_n-D)|_\infty+O(\I_n-D,2d)}{2} + \frac{|\det(\I_n+D)|_\infty+O(\I_n+D,-2d)}{2}.
\end{align*}
But over \(\ZZ_2\), \(\overline{\I_n-D} = \overline{\I_n+D}\) and \(\overline{2d} = \overline{-2d}\), hence \(O(\I_n-D,2d) =  O(\I_n+D,-2d)\). So we find the proposed formula:
\begin{equation*}
R(\varphi) = \frac{|\det(\I_n-D)|_\infty+|\det(\I_n+D)|_\infty}{2} + O(\I_n-D,2d).\qedhere
\end{equation*}
\end{proof}

\begin{proof}[Proof of \cref{prop:spectrumZZn-I}]
We will use the formula from \cref{prop:formulaRphi}. Also, recall from \cref{lem:OAaoddeven} that \(O(B,b)\) is odd (in fact, it then necessarily equals \(1\)) if and only if \(\det(B)\) is odd.

Let us first deal with the case \(n = 2\). Since \(\det(\I_2\pm D) = 1 \pm \tr (D) + \det(D)\) and \(\det(D) = \pm 1\), we have that
\begin{equation*}
\det(\I_2\pm D) \equiv \tr(D) \equiv O(\I_2-D,2d) \pmod 2.
\end{equation*}
 We now determine the value of \(R(\varphi)\):
\begin{enumerate}
	\item \(\det(D) = -1\). Then the formula becomes
	\begin{equation*}
	R(\varphi) = |\tr(D)|_\infty + O(\I_2 - D,2d).
	\end{equation*}
	Depending on the value of \(|\tr(D)|\), we have:
	\begin{enumerate}
		\item \( |\tr(D)| = 0\), then \(R(\varphi) = \infty\),
		\item \(|\tr(D)| \geq  1\), then \(R(\varphi) = |\tr(D)| + O(\I_2 - D,2d) \in 2\NN\).
	\end{enumerate}
	\item \(\det(D) = 1\). Then the formula becomes
	\begin{equation*}
	R(\varphi) = \frac{|2 - \tr (D)|_\infty+|2 + \tr (D)|_\infty}{2} + O(\I_2-D,2d).
	\end{equation*}
	Depending on the value of \(|\tr(D)|\), we have:
	\begin{enumerate}
		\item \( |\tr(D)| = 0\), then \(R(\varphi) = 2+O(\I_2-D,2d) \in 2\NN\),
		\item \(|\tr(D)| = 1\), then \(R(\varphi) = 3\),
		\item \(|\tr(D)| = 2\), then \(R(\varphi) = \infty\),
		\item \(|\tr(D)| \geq 3\), then \(R(\varphi) =  |\tr(D)|+ O(\I_2-D,2d) \in 2\NN\).
	\end{enumerate}
\end{enumerate}
So indeed \(\Spec_R(\Gamma) \subseteq 2\NN \cup \{3,\infty\}\). 
We now show that all these Reidemeister numbers can actually be attained. To obtain an even Reidemeister number, consider \(\varphi_m = \xi_{(d,D_m)}\) with 
\begin{equation*}
D_m = \begin{pmatrix}
0 & 1\\
1 & 2m
\end{pmatrix}, d= \begin{pmatrix}
1/2\\
0
\end{pmatrix},
\end{equation*}
with \(m \in \NN\), then \(|\det(\I_2-D_m)|_\infty = |\det(\I_2+D_m)|_\infty = 2|m|\) and \(O(\I_2-D_m,2d) = 0\), and hence \(R(\varphi_m) = 2|m| \). Finally, to obtain Reidemeister number \(3\), consider \(\varphi = \xi_{(0,D)}\) with
\begin{equation*}
D = \begin{pmatrix}
0 & -1\\
1 & -1
\end{pmatrix},
\end{equation*}
then \(R(\varphi) = 3\). Hence \(\Spec_R(\Gamma) = 2\NN \cup \{3,\infty\}\).

Next, consider the case \(n \geq 3\). As mentioned in the proof of \cref{prop:formulaRphi}, \(R(\varphi) \geq 2\). We show that every natural number greater than or equal to \(2\) can be attained. Consider \(\varphi = \xi_{(0,D_m)}\) with
\begin{equation*}
D_m = \begin{pmatrix}
0 	&  \cdots & \cdots & \cdots& 0 & 1\\
1 	& \ddots &&  & \vdots & 0 \\
0 & \ddots & \ddots& & \vdots & \vdots \\
\vdots & \ddots&\ddots&\ddots & \vdots &0\\
\vdots & &\ddots &\ddots& 0&m\\
0 & \cdots & \cdots& 0 & 1&m-1
\end{pmatrix},
\end{equation*}
where \(m \in \NN\). Then \(\det(\I_n - D_m) = -2m+1\), \(\det(\I_n + D_m) = (-1)^{n-1}\) and \(O(\I_n - D_m,0) = 1\), therefore \(R(\varphi) = m+1\) and thus \(\Spec_R(\Gamma) = \NN \setminus \{1\} \cup \{\infty\}\).
\end{proof}
We now determine the Reidemeister spectrum of crystallographic groups with infinite \(N_F\). For dimensions \(1\), \(2\) and \(3\), we do this for all groups; for dimension \(4\) we limit ourselves to groups where we can apply \cref{cor:directprod} or \cref{thm:averagingalmostbieb}. We will list the groups using the ``BBNWZ''-classification.

\paragraph{2/1/1/1/1} Since \(\Gamma \cong \ZZ^2\), \cref{prop:spectrumZZn} says that \(\Spec_R(\Gamma) = \NN \cup \{\infty\}\).
\paragraph{2/1/2/1/1} Since \(\Gamma \cong \langle \ZZ^2,(0,-\I_2)\rangle\), \cref{prop:spectrumZZn-I} says that \(\Spec_R(\Gamma) = 2\NN \cup \{3,\infty\}\).
\paragraph{3/1/1/1/1} Since \(\Gamma \cong \ZZ^3\), \cref{prop:spectrumZZn} says that \(\Spec_R(\Gamma) = \NN \cup \{\infty\}\).
\paragraph{3/1/2/1/1} Since \(\Gamma \cong \langle \ZZ^3,(0,-\I_3)\rangle\), \cref{prop:spectrumZZn-I} says that \(\Spec_R(\Gamma) = \NN \setminus \{1\} \cup \{\infty\}\).
\paragraph{3/2/1/1/1} Since \(\Gamma \cong \Gamma_{1/1/1/1/1} \times \Gamma_{2/1/2/1/1}\) and both factors are characteristic, \cref{cor:directprod} and \cref{tbl:roodim3} say that \(\Spec_R(\Gamma) = 4\NN \cup \{6,\infty\}\).
\paragraph{3/2/1/1/2} This is a Bieberbach group. In \cite[Proposition 4.8]{dtv18-1}, it was shown that \(\Spec_R(\Gamma) = 2\NN \cup \{\infty\}\).

\paragraph{3/2/1/2/1} This group \(\Gamma\) is given by
\begin{equation*}
\Gamma = \left\langle\ZZ^3,\left(0,\begin{pmatrix}
	1 & -1 & 0\\
	0 & -1 & 0\\
	0 & 0 & -1
\end{pmatrix}\right)\right\rangle.
\end{equation*}
We first calculate an explicit formula for the Reidemeister number of a given automorphism.

\begin{proposition}
	\label{prop:formulaRphi2}
	Let \(\Gamma = \Gamma_{3/2/1/2/1}\) and \(\varphi = \xi_{(d,D)} \in \Aut(\Gamma)\). Define
	\begin{equation*}
	\delta = \begin{cases}
	1 & \text{ if }d_3 \in \ZZ,\\
	0 & \text{ otherwise },
	\end{cases}
	\end{equation*}
	where \(d_3\) is the third coordinate of \(d\). Then
	\begin{equation*}
	R(\varphi) = \left(\frac{1}{2} \sum_{A \in F} |\det(\I_3-AD)|_\infty\right) + 4\delta.
	\end{equation*}
\end{proposition}
\begin{proof}
Let  \(\varphi = \xi_{(d,D)}\) be an automorphism of \(\Gamma\). We can calculate that \(Z(\Gamma) = \langle e_1 \rangle\) and define \(\Gamma' := \Gamma/Z(\Gamma) \cong \Gamma_{2/1/2/1/1}\). Then \(\varphi\) induces an automorphism \(\varphi' = \xi_{(d',D')}\) on \(\Gamma'\). The reader can verify, using that \(D\) commutes with any element of the holonomy group, that we may assume \(D\) and \(d\) are of the form

\begin{equation*}
D = \begin{pmatrix}
\varepsilon  & m_1 \;\;\; m_2\\
0 & D'\\
\end{pmatrix}, \;\;\;\; d = \begin{pmatrix}
0\\
d'
\end{pmatrix},
\end{equation*}
where
\begin{equation*}
D' = \begin{pmatrix}
\varepsilon + 2m_1 & 2m_2\\
m_3 & 1+2m_4\\
\end{pmatrix}, \;\;\;\; d' = \begin{pmatrix}
d_2\\
d_3
\end{pmatrix}.
\end{equation*}
Here, \(\varepsilon \in \{-1,1\}\), \(m_1, m_2, m_3, m_4, d_2 \in \ZZ\), and importantly, \(d_3 \in \frac{1}{2}\ZZ\). For \(A \in F\), let \(A'\) be the projection to the holonomy group \(F'\) of \(\Gamma'\). We have that 
\begin{equation}
\label{eq:splitepsilon}
|\det(\I_3 - AD)|_\infty = |1-\varepsilon|_\infty |\det(\I_2-A'D')|_\infty.
\end{equation}
Following  \cref{thm:det1-AD} and the first part of the proof of \cref{prop:spectrumZZn-I}, we may conclude that \(R(\varphi) = \infty\) if and only if (at least) one of the following three conditions is satisfied:
\begin{itemize}
	\item \(\varepsilon = 1\),
	\item \(\det(D') = -1\) and \(\tr(D') = 0\),
	\item  \(\det(D') = 1\) and \(|\tr(D')|=2\).
\end{itemize}
If this is the case, then the formula holds. We are left to verify the formula when none of these conditions are satisfied.

Consider a  Reidemeister class \([x]_\varphi\) of \(\Gamma\) and recall that \(Z(\Gamma) = \langle e_1 \rangle\). Then
\begin{equation*}
x = e_1^{-k}(xe_1^{2k})\varphi(e_1^{-k})^{-1}
\end{equation*}
and hence \(x \sim xe_1^{2k}\) for all \(k \in \ZZ\). So a Reidemeister class \([xZ(\Gamma)]_{\varphi'}\) of \(\Gamma'\) lifts to at most \(2\) distinct Reidemeister classes of \(\Gamma\): \([x]_\varphi\) and \([xe_1]_\varphi\).

The question that remains is: when is \(x \sim_\varphi xe_1\)? This is the case when there exists some \(z \in \Gamma\) such that
\begin{equation}
\label{eq:firstg4}
x = zxe_1\varphi(z)^{-1}.
\end{equation}
Projecting this to \(\Gamma'\) we find
\begin{equation}
\label{eq:secondg4}
xZ(\Gamma) = zx\varphi(z)^{-1}Z(\Gamma).
\end{equation}
Set \(x = ((x_1,x_2,x_3),A_x)\) and \(z = ((z_1,z_2,z_3),A_z)\). If we assume that \(A_z = \I_3\), then \eqref{eq:secondg4} is equivalent to 
\begin{equation*}
(\I_2-A_x'D')\begin{pmatrix}
z_2\\
z_3
\end{pmatrix} = 0.
\end{equation*}
But \(\det(\I_2-A_x'D') \neq 0\), hence \(z_2 = z_3 = 0\) and thus \(z = e_1^{z_1} \in Z(\Gamma)\) for some \(z_1 \in \ZZ\). But then equation \eqref{eq:firstg4} reduces to \(e_1^{2z_1+1} = 1\), which is impossible. Therefore,  \(A_z \neq \I_3\), and then equation \eqref{eq:secondg4} is a special case of equation \eqref{eq:gamma2ding}: \([xZ(\Gamma)]_{\varphi'}\) is one of the cosets of \(\I_2-A_x'D'\) such that 
\begin{equation}
\label{eq:thirdg4}
2\begin{pmatrix}
x_2\\
x_3
\end{pmatrix} + 2A'_x\begin{pmatrix}
d_2\\d_3
\end{pmatrix} = \left(\I_2-A'_xD'\right)\begin{pmatrix}
z_2\\
z_3
\end{pmatrix},
\end{equation}
i.e. a coset that forms a Reidemeister class on its own, rather than pairing up with another coset.
The \(e_2\)- and \(e_3\)-coordinates of \eqref{eq:firstg4} are equivalent to equation \eqref{eq:thirdg4}. The \(e_1\)-coordinate can be shown to be equivalent to \(z_2 = 2z_1+1\), under the assumption that \eqref{eq:thirdg4} is satisfied. But since \(z_1\) is an integer, we need that \(z_2 \in 2\ZZ+1\).

Now, let's look at the number of Reidemeister classes \([xZ(\Gamma)]_{\varphi'}\) such that \eqref{eq:thirdg4} holds. From the work we did on \(\Gamma' \cong \Gamma_{2/1/2/1/1}\), we know that we must look at the number of solutions \(O(\I_2 - D',2d')\) of the system of equations over \(\ZZ_2\) given by
\begin{equation*}
\begin{pmatrix}
0 & 0\\
\bar{m}_3 & 0
\end{pmatrix}
\begin{pmatrix}
\bar{z}_2\\
\bar{z}_3
\end{pmatrix} = \begin{pmatrix}
0 \\
\overline{2d}_3
\end{pmatrix},
\end{equation*}
and we see that \(O(\I_2 - D',2d') = 2\:O(m_3,2d_3)\), since \(\bar{z}_3\) can be chosen freely. We now have \(4\) cases:
\begin{enumerate}
	\item \(\bar{m}_3 = \bar{0}, \overline{2d}_3 = \bar{0}\). Then \(O(m_3,2d_3) = 2\) with solutions \(\bar{z}_2 = \bar{0},\;\bar{1}\); and \(\delta = 1\). 
	\item \(\bar{m}_3 = \bar{1}, \overline{2d}_3 = \bar{0}\). Then \(O(m_3,2d_3) = 1\) with solution \(\bar{z}_2 = \bar{0}\); and \(\delta = 1\).
	\item \(\bar{m}_3 = \bar{0}, \overline{2d}_3 = \bar{1}\). Then \(O(m_3,2d_3) = 0\); and \(\delta = 0\).
	\item \(\bar{m}_3 = \bar{1}, \overline{2d}_3 = \bar{1}\). Then \(O(m_3,2d_3) = 1\) with solution \(\bar{z}_2 = \bar{1}\); and \(\delta = 0\).
\end{enumerate}
Every solution \(\bar{z}_2\) of equation \(\bar{m}_3\bar{z}_2 = \overline{2d}_3\) represents \(4\) Reidemeister classes \([xZ(\Gamma)]_{\varphi'}\), since one takes all combinations of \(\bar{z}_3 \in \{0,1\}\) and \(A_x' \in \{\I_2, -\I_2\}\).
Thus, we have respectively \(8\), \(4\), \(0\) and \(4\) Reidemeister classes \([xZ(\Gamma)]_{\varphi'}\) satisfying \eqref{eq:secondg4}; of which respectively \(4\), \(0\), \(0\) and \(4\) satisfy \(z_2 \in 2\ZZ+1\). So the number of lifts to Reidemeister classes of \(\Gamma\) is respectively \(12\), \(8\), \(0\) and \(4\). This number of Reidemeister classes always equals
\begin{equation*}
2\:O(\I_2 - D',2d') + 4\delta.
\end{equation*}
On the other hand, \(\Gamma'\) has 
\begin{equation*}
\frac{|\det(\I_2-D')|+|\det(\I_2+D')|}{2} - O(\I_2-D',2d')
\end{equation*}
Reidemeister classes for which \eqref{eq:secondg4} does not hold, meaning each of these classes lift to two distinct Reidemeister classes of \(\Gamma\). Combining all the classes we obtain the formula
\begin{equation}
\label{eq:formulagamma35}
R(\varphi) = |\det(\I_2 - D')| + |\det(\I_2+D')| + 4\delta,
\end{equation}
and using \(1 - \varepsilon = 2\) in equation \eqref{eq:splitepsilon} we see this is exactly
	\begin{equation*}
R(\varphi) = \left(\frac{1}{2} \sum_{A \in F} |\det(\I_3-AD)|\right) + 4\delta = \left(\frac{1}{2} \sum_{A \in F} |\det(\I_3-AD)|_\infty\right) + 4\delta,
\end{equation*}
since none of the determinants are zero. Therefore, the proposed formula holds in all cases.
\end{proof}

\begin{proposition}
Let \(\Gamma = \Gamma_{3/2/1/2/1}\). Then \(\Spec_R(\Gamma) = 4\NN \cup \{\infty\}\).
\end{proposition}
\begin{proof}
	Let \(\varphi = \xi_{(d,D)}\) be an automorphism of \(\Gamma\) with \(R(\varphi) < \infty\). Consider formula \eqref{eq:formulagamma35} and remark that \(\tr(D') \in 2\ZZ\). Since \(\det(\I_2 \pm D') = 1 \pm \tr(D') + \det(D')\), we have
	\begin{equation*}
	|\det(\I_2 - D')| + |\det(\I_2+D')| = \begin{cases}
	4 & \text{ if } \tr(D') = 0 \text{ and } \det(D') = 1,\\
	2|\tr(D')| & \text{ otherwise },
	\end{cases}
	\end{equation*}
	so in both cases \(R(\varphi) \in 4\NN\). Now consider the automorphism \(\varphi_m = \xi_{(d,D_m)}\) given by 
\begin{equation*}
D_m = \begin{pmatrix}
-1 & m & m\\
0 & -1+2m & 2m\\
0 & 1 & 1
\end{pmatrix}, \;\;\;\; d = \begin{pmatrix}
0\\
0\\
1/2
\end{pmatrix},
\end{equation*}
where \(m \in \NN\). Since \(\det(\I_2 \pm D') = \pm 2m\) and \(\delta = 0\), we find \(R(\varphi_m) = 4|m|_\infty\). Hence \(\Spec_R(\Gamma) =  4\NN \cup \{\infty\}\).
\end{proof}

This completes the classification of the Reidemeister spectra of crystallographic groups up to dimension \(3\). We can still calculate the spectrum of some of the \(4\)-dimensional crystallographic groups with little effort:

\paragraph{4/1/1/1/1} Since \(\Gamma \cong \ZZ^4\), \cref{prop:spectrumZZn} says that \(\Spec_R(\Gamma) = \NN \cup \{\infty\}\).
\paragraph{4/1/2/1/1} Since \(\Gamma \cong \langle \ZZ^4,(0,-\I_4)\rangle\), \cref{prop:spectrumZZn-I} says that \(\Spec_R(\Gamma) = \NN \setminus \{1\} \cup \{\infty\}\).
\paragraph{4/2/2/1/1} Since \(\Gamma \cong \Gamma_{1/1/1/1/1} \times \Gamma_{3/1/2/1/1}\) and both factors are characteristic, \cref{cor:directprod} says that \(\Spec_R(\Gamma) = 2\NN \setminus \{2\} \cup \{\infty\}\).
\paragraph{4/2/2/1/2} This is a Bieberbach group. In \cite[Proposition 4.8]{dtv18-1}, it was shown that \(\Spec_R(\Gamma) = 2\NN \cup \{\infty\}\).
\paragraph{4/3/1/1/1} Since \(\Gamma \cong \Gamma_{2/1/1/1/1} \times \Gamma_{2/1/2/1/1}\) and both factors are characteristic, \cref{cor:directprod} says that \(\Spec_R(\Gamma) = 2\NN \cup 3\NN \cup \{\infty\}\).
\paragraph{4/3/1/1/2} This is a Bieberbach group. In \cite[Proposition 5.10, case \(\delta = 0\)]{dtv18-1}, it was shown that \(\Spec_R(\Gamma) = 2\NN \cup \{\infty\}\).
\paragraph{4/3/1/2/2} This is a Bieberbach group. In \cite[Proposition 5.10, case \(\delta = 1\)]{dtv18-1}, it was shown that \(\Spec_R(\Gamma) = 4\NN \cup \{\infty\}\).
\paragraph{4/8/1/1/2} This is a Bieberbach group. In \cite[Proposition 5.12, case \(\delta = 1\)]{dtv18-1}, it was shown that \(\Spec_R(\Gamma) = 6\NN \cup \{\infty\}\).
\paragraph{4/8/1/2/1} Since \(\Gamma \cong \Gamma_{2/1/1/1/1} \times \Gamma_{2/4/1/1/1}\) and both factors are characteristic, \cref{cor:directprod} says that \(\Spec_R(\Gamma) = 4\NN \cup \{\infty\}\).
\paragraph{4/8/1/2/2} This is a Bieberbach group. In \cite[Proposition 5.10, case \(\delta = 0\)]{dtv18-1}, it was shown that \(\Spec_R(\Gamma) = 6\NN \cup \{\infty\}\).
\paragraph{4/9/2/1/1} Since \(\Gamma \cong \Gamma_{2/1/2/1/1} \times \Gamma_{2/4/1/1/1}\) and both factors are characteristic, \cref{cor:directprod} says that \(\Spec_R(\Gamma) = 8\NN \cup \{12,\infty\}\).

\section{Conclusion}
The first algorithm given in this paper calculates whether for a given \(n\)-dimensional crystallographic group \(\Gamma\) with holonomy group \(F\) and a given \(D \in N_F\), there exists a \(d \in \RR^n\) such that the map \(\gamma \mapsto (d,D)\gamma(d,D)^{-1}\) is an automorphism of \(\Gamma\). 

We then presented two main algorithms for crystallographic groups whose holonomy group \(F\) has finite normaliser in \(\GL_{n}(\ZZ)\). \Cref{alg:mainalg1} determines if a given group admits the \(R_\infty\)-property; \cref{alg:Reidnr} determines the Reidemeister spectrum. The former main algorithm was applied to the crystallographic groups up to dimension \(6\), leading to the results in \cref{tbl:alg2results}. The latter was applied to crystallographic groups up to dimension \(4\). We complemented the output of both algorithms with results on crystallographic groups whose holonomy group has infinite normaliser, to obtain the data found in \cref{tbl:roodim3} and \cref{tbl:roodim4}. 

In particular we have obtained complete classifications of the crystallographic groups without the \(R_\infty\)-property up to dimension \(4\), of the Reidemeister spectra of crystallographic groups up to dimension \(3\), and of the Reidemeister spectra of Bieberbach groups up to dimension \(4\).

\begin{table}[h]
	\centering
	\begin{tabular}{l|l|l|l|l}
		BBNWZ & IT & CARAT & \(|N_F|\) & \(\Spec_R(\Gamma)\setminus \{\infty\}\)\\
		\hline
		\hline
		&&&& \\[\dimexpr-\normalbaselineskip+2pt]
		\(1/1/1/1/1*\) & 1/1	 &min.1-1.1-0& \(2\) &  \(\{2\}\)  \\
		\hline
		&&&& \\[\dimexpr-\normalbaselineskip+2pt]
		\(2/1/1/1/1*\) 	&2/1  	 & min.2-1.1-0& \(\infty\)& \(\NN\)  \\
		\(2/1/2/1/1\) 	&2/2 	 & group.1-1.1-0& \(\infty\)& \(2\NN \cup \{3\}\) \\
		\(2/4/1/1/1\)		&2/13	 &min.5-1.1-0& \(12\) &  \(\{4\}\)\\
		\hline
		&&&& \\[\dimexpr-\normalbaselineskip+2pt]
		\(3/1/1/1/1*\) 	&3/1 	 &min.6-1.1-0& \(\infty\) &  \(\NN\)\\
		\(3/1/2/1/1\) 	&3/2 	 &group.5-1.1-0& \(\infty\)&  \(\NN\setminus \{1\}\)\\
		\(3/2/1/1/1\)		&3/3  	& min.7-1.1-0& \(\infty\)&  \(4\NN \cup \{6\}\)\\
		\(3/2/1/1/2*\)	&3/4  	 &min.7-1.1-1& \(\infty\)&  \(2\NN\)\\
		\(3/2/1/2/1\)		&3/5  	 &min.7-1.2-0& \(\infty\)&  \(4\NN\)\\
		\(3/3/1/1/1\)		&3/16 	 &min.10-1.1-0& \(48\)&  \(\{2\}\)  \\
		\(3/3/1/1/4*\)	&3/19 &min.10-1.1-3& \(48\)&  \(\{2\}\)  \\
		\(3/3/1/3/1\)		&3/22 &min.10-1.3-0& \(48\)&  \(\{2\}\) \\
		\(3/3/1/4/1\)		&3/23 &min.10-1.4-0& \(48\)&  \(\{2\}\) \\
		\(3/3/1/4/2\)		&3/24&min.10-1.4-1& \(48\) &  \(\{2\}\)\\
		\(3/5/1/1/1\)		&3/146  &min.13-1.2-0& \(12\)& \(\{8\}\)  \\
		\(3/5/1/2/1\)		&3/143   &min.13-1.1-0& \(24\)& \(\{8\}\)
	\end{tabular}
\captionsetup{justification=centering}
	\caption{\(1\)-, \(2\)- and \(3\)-dimensional crystallographic groups without \(R_\infty\)\\(Note that we have omitted the value \(\infty\) from the Reidemeister spectra)}
	\label{tbl:roodim3}
\end{table}

\begin{table}[!ht]
	\centering
	\tiny 
	\begin{tabular}{l|l|c|c}
		BBNWZ & CARAT &\(|N_F|\) & \(\Spec_R(\Gamma)\!\setminus\! \{\infty\}\)\\
		\hline
		\hline
		&&& \\[\dimexpr-\normalbaselineskip+2pt]
		\(4/1/1/1/1\)*	& min.15-1.1-0&\(\infty\)	& \(\NN\)\\
		\(4/1/2/1/1\)	& group.26-1.1-0&\(\infty\)	& \(\NN \setminus \{1\}\)\\
		\(4/2/2/1/1\)	& min.17-1.1-0&\(\infty\)	& \(2\NN \setminus \{2\}\)\\
		\(4/2/2/1/2\)*	& min.17-1.1-1&\(\infty\)	& \(2\NN\)\\
		\(4/2/2/2/1\)	& min.17-1.2-0&\(\infty\)	& \\
		
		\(4/3/1/1/1\)	& min.18-1.1-0&\(\infty\)	& \(2\NN \cup 3\NN\)\\
		\(4/3/1/1/2\)*	& min.18-1.1-1&\(\infty\)	& \(2\NN\)\\
		\(4/3/1/2/1\)	& min.18-1.2-0&\(\infty\)	& \\
		\(4/3/1/2/2\)*	& min.18-1.2-1&\(\infty\)	& \(4\NN\)\\
		\(4/3/1/3/1\)	& min.18-1.3-0&\(\infty\)	& \\
		
		\(4/3/2/1/1\)	& group.28-1.1-0&\(\infty\)	& \\
		\(4/3/2/1/2\)	& group.28-1.1-1&\(\infty\)	& \\
		\(4/3/2/1/3\) 	& group.28-1.1-2&\(\infty\)	& \\
		\(4/3/2/2/1\) 	& group.28-1.2-0&\(\infty\)	& \\
		\(4/3/2/2/2\)	& group.28-1.2-1&\(\infty\)	& \\
		
		\(4/3/2/2/3\)	& group.28-1.2-2&\(\infty\)	& \\
		\(4/3/2/3/1\)	& group.28-1.3-0&\(\infty\)	& \\
		\(4/5/1/1/1\)	& group.52-1.13-0&\(48\)		& \(\{4\}\)\\
		\(4/5/1/2/1\)	& group.52-1.1-0&\(96\)		& \(\{4\}\)\\
		\(4/5/1/2/9\)*	& group.52-1.1-6&\(96\)		& \(\{4\}\)\\
		
		\(4/5/1/5/1\)	& group.52-1.7-0&\(96\)		& \(\{4\}\)\\
		\(4/5/1/5/2\)	& group.52-1.7-1&\(96\)		& \(\{4\}\)\\
		\(4/5/1/7/1\)	& group.52-1.12-0&\(96\)		& \(\{4\}\)\\
		\(4/5/1/7/4\)*	& group.52-1.12-3&\(96\)		& \(\{4\}\)\\
		\(4/5/1/9/1\)	& group.52-1.3-0&\(96\)		& \(\{4\}\)\\
		
		\(4/5/1/13/1\)	& group.52-1.6-0&\(96\)		& \(\{4\}\)\\
		\(4/5/2/2/1\)	& group.80-1.1-0&\(384\)		& \(\{2,4\}\)\\
		\(4/5/2/2/16\)	& group.80-1.1-5&\(384\)		& \(\{2,4\}\)\\
		\(4/5/2/2/18\)	& group.80-1.1-18&\(384\)		& \(\{2,4\}\)\\
		\(4/5/2/5/1\)	& group.80-1.8-0&\(384\)		& \(\{4\}\)\\
		
		\(4/5/2/5/3\)	& group.80-1.8-4&\(384\)		& \(\{2\}\)\\
		\(4/5/2/5/5\)	& group.80-1.8-2&\(384\)		& \(\{4\}\)\\
		\(4/5/2/5/6\)	& group.80-1.8-5&\(384\)		& \(\{2\}\)\\
		\(4/5/2/6/1\)	& group.80-1.6-0&\(128\)		& \(\{4\}\)\\
		\(4/5/2/6/3\)	& group.80-1.6-2&\(128\)		& \(\{4\}\)\\
		
		\(4/5/2/9/1\)	& group.80-1.4-0&\(384\)		& \(\{4\}\)\\
		\(4/5/2/9/3\)	& group.80-1.4-2&\(384\)		& \(\{4\}\)\\
		\(4/8/1/1/1\)	& group.179-1.2-0&\(\infty\)	& \\
		\(4/8/1/1/2\)*	& group.179-1.2-1&\(\infty\)	& \(6\NN\)\\
		\(4/8/1/2/1\)	& group.179-1.1-0&\(\infty\)	& \(4\NN\)\\
		
		\(4/8/1/2/2\)*	& group.179-1.1-1&\(\infty\)	& \(6\NN\)\\
		\(4/9/2/1/1\)	& group.182-1.1-0&\(\infty\)	& \(8\NN \cup \{12\}\)\\
		\(4/10/1/1/1\)	& min.36-1.1-0&\(\infty\)	& \\
		\(4/11/1/1/1\)	& group.170-1.1-0&\(\infty\)	& \\
		\(4/11/2/1/1\)	& group.171-1.1-0&\(\infty\)	& \\
		
		\(4/16/1/1/1\)	& group.96-1.1-0&\(\infty\)	& \\
		
	\end{tabular}
	\begin{tabular}{l|l|c|c}
		BBNWZ & CARAT &\(|N_F|\) & \(\Spec_R(\Gamma)\!\setminus\! \{\infty\}\)\\
		\hline
		\hline
		&&& \\[\dimexpr-\normalbaselineskip+2pt]
		\(4/16/1/1/2\)	& group.96-1.1-1&\(\infty\)	&\\
		\(4/16/1/2/1\)	& group.96-2.1-0&\(\infty\)	&\\
		\(4/16/1/2/2\)	& group.96-2.1-1&\(\infty\)	&\\
		\(4/16/1/2/3\)	& group.96-2.1-2&\(\infty\)	&\\
		\(4/16/1/3/1\)	& group.96-3.1-0&\(\infty\)	&\\
		
		\(4/17/1/1/1\)	& group.173-2.1-0&\(\infty\)	&\\
		\(4/17/1/2/1\)	& group.173-3.1-0&\(\infty\)	&\\
		\(4/17/1/3/1\)	& group.173-1.1-0&\(\infty\)	&\\
		\(4/17/2/1/1\)	& group.172-2.1-0&\(\infty\)	&\\
		\(4/17/2/2/1\)	& group.172-1.1-0&\(\infty\)	&\\
		
		\(4/18/1/2/1\)	& group.169-1.1-0&\(128\)		&\(\{4,8\}\)\\
		\(4/18/1/2/3\)	& group.169-1.1-2&\(128\)		&\(\{4,8\}\)\\
		\(4/18/1/3/1\)	& group.169-1.2-0&\(128\)		&\(\{4,8\}\)\\
		\(4/18/1/3/2\)	& group.169-1.2-1&\(128\)		&\(\{4,8\}\)\\
		\(4/18/4/2/1\)	& group.163-1.1-0&\(128\)		&\(\{4,8\}\)\\
		
		\(4/18/4/2/3\)	& group.163-1.1-6&\(128\)		&\(\{4,8\}\)\\
		\(4/18/4/2/6\)	& group.163-1.1-4&\(128\)		&\(\{4,8\}\)\\
		\(4/18/4/5/1\)	& group.163-1.2-0&\(128\)		&\(\{4,8\}\)\\
		\(4/18/4/5/3\)	& group.163-1.2-2&\(128\)		&\(\{4,8\}\)\\
		\(4/18/4/5/5\)	& group.163-1.2-7&\(128\)		&\(\{4,8\}\)\\
		
		\(4/18/4/5/6\)	& group.163-1.2-6&\(128\)		&\(\{4,8\}\)\\
		\(4/21/2/2/1\)	& group.37-1.1-0&\(144\)		&\(\{3\}\)\\
		\(4/22/1/1/1\)	& min.32-1.2-0 &\(144\)		&\(\{16\}\)\\
		\(4/22/1/2/1\)	& min.32-1.1-0&\(288\)		&\(\{4,16\}\)\\
		\(4/22/2/2/1\)	& group.40-1.1-0&\(288\)		&\(\{8\}\)\\
		
		\(4/22/5/3/1\)	& group.44-3.1-0&\(288\)		&\(\{6\}\)\\
		\(4/22/5/4/1\)	& group.44-1.1-0&\(288\)		&\(\{6\}\)\\
		\(4/22/7/2/1\)	& min.28-1.1-0&\(288\)		&\(\{12\}\)\\
		\(4/26/1/1/1\)	& group.109-1.1-0&\(\infty\)	&\\
		\(4/26/2/1/1\)	& max.6-1.1-0&\(\infty\)	&\\
		
		\(4/26/2/1/2\)	& max.6-1.1-1&\(\infty\)	&\\
		\(4/27/1/1/1\)	& group.144-1.1-0&\(\infty\)	&\\
		\(4/27/2/1/1\)	& group.141-1.1-0&\(\infty\)	&\\
		\(4/27/3/1/1\)	& group.142-2.1-0&\(\infty\)	&\\
		\(4/27/3/2/1\)	& group.142-1.1-0&\(\infty\)	&\\
		
		\(4/27/4/1/1\)	& group.143-1.1-0&\(\infty\)	&\\
		\(4/28/1/1/1\)	& min.43-1.1-0&\(\infty\)	&\\
		\(4/28/2/1/1\)	& min.44-1.1-0&\(\infty\)	&\\
		\(4/32/1/2/1\)	& group.103-1.1-0&\(576\)		&\(\{2,6\}\)\\
		\(4/32/1/2/2\)	& group.103-1.1-1&\(576\)		&\(\{2,6\}\)\\
		
		\(4/32/4/2/1\)	& group.78-1.1-0&\(192\)		&\(\{2,6\}\)\\
		\(4/32/4/2/3\)	& group.78-1.1-2&\(192\)		&\(\{2,6\}\)\\
		\(4/32/4/2/6\)	& group.78-1.1-4&\(192\)		&\(\{2,6\}\)\\
		\(4/32/10/2/1\)	& min.38-1.1-0&\(1152\)		&\(\{6\}\)\\
		\(4/32/10/2/7\)	& min.38-1.1-4&\(1152\)		&\(\{6\}\)\\
		\multicolumn{3}{c}{}
	\end{tabular}
	\captionsetup{justification=centering}
	\caption{\(4\)-dimensional crystallographic groups without \(R_\infty\)\\(Note that we have omitted the value \(\infty\) from the Reidemeister spectra)}
	\label{tbl:roodim4}
\end{table}

\paragraph{Acknowledgement} The authors would like to thank prof. Bettina Eick for suggesting the use of GAP in studying crystallographic groups, prof. Franz Gähler for answering some questions about CARAT, and the referee for their detailed remarks and suggestions for the paper. 

\printbibliography

\noindent Karel Dekimpe, Sam Tertooy\\
KU Leuven Campus Kulak Kortrijk\\
E.\ Sabbelaan 53\\
8500 Kortrijk\\
Belgium\\[2mm]
Karel.Dekimpe\@@kuleuven.be\\
Sam.Tertooy\@@kuleuven.be\\
\phantom{Phantom} \\
Tom Kaiser\\
Universit\'e de Neuch\^{a}tel\\
Institut de Math\'ematiques\\
Rue Emile--Argand 11\\
2000 Neuch\^{a}tel\\
Switzerland\\[2mm]
Tom.Kaiser\@@unine.ch
\end{document}